\documentclass[12pt]{article}
\usepackage{amsmath,amssymb,amsbsy,amsfonts,amsthm,latexsym,amsopn,amstext,
                                               amsxtra,euscript,amscd}
\begin{document}

\newfont{\teneufm}{eufm10}
\newfont{\seveneufm}{eufm7}
\newfont{\fiveeufm}{eufm5}
%
%
\newfam\eufmfam
                \textfont\eufmfam=\teneufm \scriptfont\eufmfam=\seveneufm
                \scriptscriptfont\eufmfam=\fiveeufm
%
%
\def\frak#1{{\fam\eufmfam\relax#1}}
%


\def\bbbr{{\rm I\!R}} 
\def\bbbm{{\rm I\!M}}
\def\bbbn{{\rm I\!N}} 
\def\bbbf{{\rm I\!F}}
\def\bbbh{{\rm I\!H}}
\def\bbbk{{\rm I\!K}}
\def\bbbp{{\rm I\!P}}
\def\bbbone{{\mathchoice {\rm 1\mskip-4mu l} {\rm 1\mskip-4mu l}
{\rm 1\mskip-4.5mu l} {\rm 1\mskip-5mu l}}}
\def\bbbc{{\mathchoice {\setbox0=\hbox{$\displaystyle\rm C$}\hbox{\hbox
to0pt{\kern0.4\wd0\vrule height0.9\ht0\hss}\box0}}
{\setbox0=\hbox{$\textstyle\rm C$}\hbox{\hbox
to0pt{\kern0.4\wd0\vrule height0.9\ht0\hss}\box0}}
{\setbox0=\hbox{$\scriptstyle\rm C$}\hbox{\hbox
to0pt{\kern0.4\wd0\vrule height0.9\ht0\hss}\box0}}
{\setbox0=\hbox{$\scriptscriptstyle\rm C$}\hbox{\hbox
to0pt{\kern0.4\wd0\vrule height0.9\ht0\hss}\box0}}}}
\def\bbbq{{\mathchoice {\setbox0=\hbox{$\displaystyle\rm
Q$}\hbox{\raise 0.15\ht0\hbox to0pt{\kern0.4\wd0\vrule
height0.8\ht0\hss}\box0}} {\setbox0=\hbox{$\textstyle\rm
Q$}\hbox{\raise 0.15\ht0\hbox to0pt{\kern0.4\wd0\vrule
height0.8\ht0\hss}\box0}} {\setbox0=\hbox{$\scriptstyle\rm
Q$}\hbox{\raise 0.15\ht0\hbox to0pt{\kern0.4\wd0\vrule
height0.7\ht0\hss}\box0}} {\setbox0=\hbox{$\scriptscriptstyle\rm
Q$}\hbox{\raise 0.15\ht0\hbox to0pt{\kern0.4\wd0\vrule
height0.7\ht0\hss}\box0}}}}
\def\bbbt{{\mathchoice {\setbox0=\hbox{$\displaystyle\rm
T$}\hbox{\hbox to0pt{\kern0.3\wd0\vrule height0.9\ht0\hss}\box0}}
{\setbox0=\hbox{$\textstyle\rm T$}\hbox{\hbox
to0pt{\kern0.3\wd0\vrule height0.9\ht0\hss}\box0}}
{\setbox0=\hbox{$\scriptstyle\rm T$}\hbox{\hbox
to0pt{\kern0.3\wd0\vrule height0.9\ht0\hss}\box0}}
{\setbox0=\hbox{$\scriptscriptstyle\rm T$}\hbox{\hbox
to0pt{\kern0.3\wd0\vrule height0.9\ht0\hss}\box0}}}}
\def\bbbs{{\mathchoice
{\setbox0=\hbox{$\displaystyle     \rm S$}\hbox{\raise0.5\ht0\hbox
to0pt{\kern0.35\wd0\vrule height0.45\ht0\hss}\hbox
to0pt{\kern0.55\wd0\vrule height0.5\ht0\hss}\box0}}
{\setbox0=\hbox{$\textstyle        \rm S$}\hbox{\raise0.5\ht0\hbox
to0pt{\kern0.35\wd0\vrule height0.45\ht0\hss}\hbox
to0pt{\kern0.55\wd0\vrule height0.5\ht0\hss}\box0}}
{\setbox0=\hbox{$\scriptstyle      \rm S$}\hbox{\raise0.5\ht0\hbox
to0pt{\kern0.35\wd0\vrule height0.45\ht0\hss}\raise0.05\ht0\hbox
to0pt{\kern0.5\wd0\vrule height0.45\ht0\hss}\box0}}
{\setbox0=\hbox{$\scriptscriptstyle\rm S$}\hbox{\raise0.5\ht0\hbox
to0pt{\kern0.4\wd0\vrule height0.45\ht0\hss}\raise0.05\ht0\hbox
to0pt{\kern0.55\wd0\vrule height0.45\ht0\hss}\box0}}}}
\def\bbbz{{\mathchoice {\hbox{$\sf\textstyle Z\kern-0.4em Z$}}
{\hbox{$\sf\textstyle Z\kern-0.4em Z$}} {\hbox{$\sf\scriptstyle
Z\kern-0.3em Z$}} {\hbox{$\sf\scriptscriptstyle Z\kern-0.2em
Z$}}}}
\def\ts{\thinspace}

\newtheorem{theorem}{Theorem}
\newtheorem{lemma}[theorem]{Lemma}
\newtheorem{claim}[theorem]{Claim}
\newtheorem{cor}[theorem]{Corollary}
\newtheorem{prop}[theorem]{Proposition}
\newtheorem{definition}{Definition}
\newtheorem{question}[theorem]{Open Question}

\def\squareforqed{\hbox{\rlap{$\sqcap$}$\sqcup$}}
\def\qed{\ifmmode\squareforqed\else{\unskip\nobreak\hfil
\penalty50\hskip1em\null\nobreak\hfil\squareforqed
\parfillskip=0pt\finalhyphendemerits=0\endgraf}\fi}

\def\cA{{\mathcal A}}
\def\cB{{\mathcal B}}
\def\cC{{\mathcal C}}
\def\cD{{\mathcal D}}
\def\cE{{\mathcal E}}
\def\cF{{\mathcal F}}
\def\cG{{\mathcal G}}
\def\cH{{\mathcal H}}
\def\cI{{\mathcal I}}
\def\cJ{{\mathcal J}}
\def\cK{{\mathcal K}}
\def\cL{{\mathcal L}}
\def\cM{{\mathcal M}}
\def\cN{{\mathcal N}}
\def\cO{{\mathcal O}}
\def\cP{{\mathcal P}}
\def\cQ{{\mathcal Q}}
\def\cR{{\mathcal R}}
\def\cS{{\mathcal S}}
\def\cT{{\mathcal T}}
\def\cU{{\mathcal U}}
\def\cV{{\mathcal V}}
\def\cW{{\mathcal W}}
\def\cX{{\mathcal X}}
\def\cY{{\mathcal Y}}
\def\cZ{{\mathcal Z}}




\def\sssum{\mathop{\sum\!\sum\!\sum}}
\def\ssum{\mathop{\sum\!\sum}}

\newcommand{\ignore}[1]{}

\def\vec#1{\mathbf{#1}}

\def\e{\mathbf{e}}



\def\AA{\mathbb{A}}
\def\BB{\mathbf{B}}

\def\GL{\mathrm{GL}}

\hyphenation{re-pub-lished}

\def\rank{{\mathrm{rk}\,}}
\def\ad{{\mathrm ad}}

\def \C{\mathbb{C}}
\def \F{\mathbb{F}}
\def \K{\mathbb{K}}
\def \Z{\mathbb{Z}}
\def \R{\mathbb{R}}
\def \Q{\mathbb{Q}}
\def \N{\mathbb{N}}

\def \nd{{\, | \hspace{-1.5 mm}/\,}}

\def\Zn{\Z_n}

\def\Fp{\F_p}
\def\Fq{\F_q}
\def \fp{\Fp^*}
\def\\{\cr}
\def\({\left(}
\def\){\right)}
\def\fl#1{\left\lfloor#1\right\rfloor}
\def\rf#1{\left\lceil#1\right\rceil}

\def\Spec#1{\mbox{\rm {Spec}}\,#1}
\def\invp#1{\mbox{\rm {inv}}_p\,#1}
\def\Mp{\cM_{m,n}(\F_p)}
\def\MpH{\cM_{m,n}(H;\F_p)}
\def\Mpn{\cM_{n}(\F_p)}
\def\MpHn{\cM_{n}(H;\F_p)}
\def\Rp{R_{m,n}(k;\F_p)}
\def\rp{r_{m,n}(k;\F_p)}
\def\RpH{R_{m,n}(k,H;\F_p)}
\def\FpH{F_{n,\vt}(H;\F_p)}
\def\UpH{U_{n}(H;\F_p)}
\def\MZ{\cM_n(\Z)}
\def\vt{\vec{t}}

\def\MpD{\widetilde\cM_{m,n}(D;\F_p)}
\def\RpD{\widetilde R_{m,n}(k,D;\F_p)}
\def\FpD{\widetilde F_{n,\vt}(D;\F_p)}
\def\UpD{\widetilde U_{n}(D;\F_p)}

\def\SL{\mathrm{SL}}

\def\Ln#1{\mbox{\rm {Ln}}\,#1}

\def\epp{\mbox{\bf{e}}_{p-1}}
\def\ep{\mbox{\bf{e}}_p}
\def\em{\mbox{\bf{e}}_{m}}
\def\ed{\mbox{\bf{e}}_{d}}

\def\ii {\iota}

\def\wt#1{\mbox{\rm {wt}}\,#1}

\def\GR#1{{ \langle #1 \rangle_n }}

\def\ab{\{\pm a,\pm b\}}
\def\cd{\{\pm c,\pm d\}}

\def\Bt {\mbox{\rm {Bt}}}

\def\Res#1{\mbox{\rm {Res}}\,#1}

\def\Tr#1{\mbox{\rm {Tr}}\,#1}

\setlength{\textheight}{43pc}
\setlength{\textwidth}{28pc}

\title{On The Solvability of Bilinear Equations in Finite 
Fields}

 \author{ 
 {\sc Igor E.~Shparlinski}\\
 {Department of Computing}\\
 {Macquarie University, Sydney, NSW 2109, Australia} \\
{\tt igor@ics.mq.edu.au} }

\maketitle

\begin{abstract} We consider the equation 
$$
ab + cd = \lambda, \qquad a\in\cA,\ b \in \cB, \ c\in \cC, \ d \in \cD , 
$$
over a finite field $\F_q$ of $q$ elements, with variables from arbitrary sets
$\cA,\cB, \cC,  \cD \subseteq  \F_q$.
The question of solvability of such and more general equations 
has recently
been  considered by D.~Hart and A.~Iosevich, who, in particular,
proved that  
$$
\# \cA  \# \cB \# \cC \#  \cD \ge C   q^3
$$
for some absolute constant $C>0$, 
 then above equation has a solution for any $\lambda \in \F_q^*$. 
Here we show that using 
bounds of multiplicative character sums allows us to extend the 
class of sets which satisfy this property. 
\end{abstract}


\newpage

\section{Introduction} 

\subsection{Background}

Let $\F_q$ denote the final field of $q$ elements. 

Using exponential sums,  D. Hart and A. Iosevich~\cite{HaIo} 
have shown that for any $2n$ sets $\cA_i,
\cB_i\subseteq  \F_q$, $i=1, \ldots,n $, with 
\begin{equation}
\label{eq:HaIo n}
\prod_{i=1}^n \# \cA_i  \# \cB_i > C q^{n+1}
\end{equation} 
 for some absolute constant $C > 0$, 
the equation
\begin{equation}
\label{eq:sum ab}
\sum_{i=1}^n a_ib_i= \lambda, \qquad a_i\in\cA_i,\ b_i \in \cB_i, \ i=1, \ldots,n , 
\end{equation} 
has a solution for any $\lambda \in \F_q^*$ (although the proof is given  
only in the case  of $\cA_1 = \cB_1= \ldots = \cA_n = \cB_n$, 
the method and results immediately
extend to the general case, see~\cite[Remark~1.3]{HaIo}). These results have been put 
in a more general context in a recent work of  D. Hart A. Iosevich, D. Koh and
M. Rudnev~\cite{HaIoKoRu}.

In particular, for $n=2$, one can easily derive from the proof 
of~\cite[Theorem~1.1]{HaIo} that 
the equation
\begin{equation}
\label{eq:ab+cd}
ab + cd = \lambda, \qquad a\in\cA,\ b \in \cB, \ c\in \cC, \ d \in \cD, 
\end{equation} 
has
\begin{equation}
\label{eq:HaIo N}
N(\cA, \cB, \cC, \cD;\lambda) =   \frac{\# \cA  \# \cB \# \cC \#  \cD }{q-1}
 + O\(\(q\# \cA  \# \cB \#\cC \#  \cD\)^{1/2}\)
\end{equation} 
solutions for any $\lambda \in \F_q^*$.

In particular, we see from~\eqref{eq:HaIo n} that  for any $\lambda \in \F_q^*$,
the equation~\eqref{eq:ab+cd} has a solution  
for any four  sets $\cA, \cB, \cC, \cD \subseteq  \F_q^*$
with 
\begin{equation}
\label{eq:HaIo 2}
\# \cA  \# \cB \# \cC \#  \cD \ge C   q^3. 
\end{equation} 
for some absolute constant $C>0$.

Furthermore, for the number $T(\cA, \cB, \cC, \cD)$ 
of the solutions of the similar  equation 
\begin{equation}
\label{eq:a+b+cd}
a + b = cd, \qquad a \in \cA, \ b\in  \cB,\ c \in \cC,\ d \in \cD,
\end{equation} 
 A.~S{\'a}rk{\"o}zy~\cite{Sark} has  given an analogous bound
\begin{equation}
\label{eq:Sark T}
T(\cA, \cB, \cC, \cD)  =   \frac{\# \cA  \# \cB \# \cC \#  \cD }{q-1}
 + O\(\(q\# \cA  \# \cB \#\cC \#  \cD\)^{1/2}\)
\end{equation}
(it is shown in~\cite{Sark} only for prime fields
but the proof and result extend  to the general case at the 
cost of only typographical changes). 
Thus~\eqref{eq:a+b+cd} has a solution 
under the condition~\eqref{eq:HaIo 2} as well. 

D. Hart and A. Iosevich~\cite{HaIo} have also considered the set 
$\cE(\cA, \cB, \cC, \cD)$
of $\lambda \in \F_q$ for which~\eqref{eq:ab+cd} has no solution. 
Although only the case of $\cA = \cB = \cC = \cD$ 
has been worked out in~\cite{HaIo}, the same approach seems to give  the
bound
$$
\# \cE(\cA, \cB, \cC, \cD) = O\(\frac{q^3}{ \# A \# \cB \# \cC }\)
$$
(certainly if $\cD$ is not of the
smallest cardinality among $\cA, \cB, \cC, \cD$,
one can alter this bound in the obvious way).
Thus only three sets out four are relevant.
This implies that instead of~\eqref{eq:ab+cd} one can consider 
the equation with only three variables
\begin{equation}
\label{eq:f+gh}
f + gh = \lambda, \qquad f \in \cF, \ g\in  \cG,\ h \in \cH,  
\end{equation} 
for some sets $ \cF,\cG, \cH \subseteq \F_q$. 
If   $\cE(\cF,\cG, \cH)$ is
the set  of $\lambda \in \F_q$ for which~\eqref{eq:f+gh} has no solution
then  we have
\begin{equation}
\label{eq:EABC}
\#  \cE(\cF,\cG, \cH) = O\(\frac{q^3}{ \# \cF \# \cG \# \cH }\).
\end{equation}

We now remark that the result of A.~S{\'a}rk{\"o}zy~\cite{Sark} also
implies~\eqref{eq:EABC}. 
Indeed, we see that 
$$
T\(-\cE(\cF,\cG, \cH),    \cF,   \cG, \cH\) = 0
$$
and we derive from~\eqref{eq:Sark T} that
$$
\frac{\#  \cE(\cF,\cG, \cH)  \# \cF \# \cG \# \cH  }{q-1} =
  O\(\(q \# \cE(\cF,\cG, \cH) \# \cF \# \cG \# \cH  \)^{1/2}\),
$$
which leads to~\eqref{eq:EABC}. 

Furthermore, in a recent work of M.~Garaev~\cite{Gar2}
considering the equation~\eqref{eq:a+b+cd} (for some 
special sets $\cA, \cB, \cC, \cD$ has been the main tool
in obtaining new results on the sum-product problem in 
finite fields (see also~\cite{BGK,BKT,Gar1,HaIo,HaIoSo,KatzShen1,
KatzShen2} for the background). In fact, one can shorten  
the proof of~\cite[Theorem~1]{Gar2} by a direct appeal
to~\eqref{eq:Sark T}. 

\subsection{Our Results}

Here we show that using multiplicative character sums
one can extend the region of possible values 
for  $\# \cA,  \# \cB,  \# \cC, \#  \cD$
which guarantee the solvability of  
the equations~\eqref{eq:ab+cd} and~\eqref{eq:a+b+cd}. 
We show that  for any fixed 
$\varepsilon>0$, there exists $\delta > 0$ such that if   
\begin{equation}
\label{eq:A and B}
\# \cA  \ge q^{1/2 + \varepsilon} \qquad \text{and} 
\qquad  \# \cB \ge q^{\varepsilon},  
\end{equation} 
as well as
\begin{equation}
\label{eq:C and D}
\# \cC \# \cD \ge  q^{2-\delta}, 
\end{equation} 
then the equations~\eqref{eq:ab+cd} and~\eqref{eq:a+b+cd}
have a solution for any $\lambda \in \F_q^*$ (provided that $q$ 
is large enough). 

More precisely, we obtain the following asymptotic formulas:

\begin{theorem}
\label{thm:2 prod}
For any fixed 
$\varepsilon>0$, there exists $\delta > 0$ such that
for  any $\lambda \in \F_q^*$  and any 
sets $\cA, \cB, \cC, \cD \subseteq  \F_q$ which
satisfy ~\eqref{eq:A and B} and~\eqref{eq:C and D},  the 
equations~\eqref{eq:ab+cd} and~\eqref{eq:a+b+cd}
have 
$$
N(\cA, \cB, \cC, \cD;\lambda) =  \frac{\# \cA  \# \cB \# \cC \#  \cD}{q} \(1  +  O(q^{-\delta})\)
$$
and
$$
T(\cA, \cB, \cC, \cD) =  \frac{\# \cA  \# \cB \# \cC \#  \cD}{q} \(1  +  O(q^{-\delta})\)
$$
  solutions, respectively.  
\end{theorem}

The same argument which we have used to derive~\eqref{eq:EABC}
from~\eqref{eq:Sark T}, combined with 
our bound on  $T(\cA, \cB, \cC, \cD)$, 
leads to the following estimates:

\begin{cor}
\label{cor:EFGH1}
For any fixed 
$\varepsilon>0$, there exists $\delta > 0$ such that for  any  
sets $ \cF,\cG, \cH \subseteq  \F_q$ which
satisfy 
$$
  \# \cF  \ge q^{  \varepsilon} \qquad \text{and}  \qquad 
\# \cG \# \cH \ge  q^{2-\delta}
$$
the 
equation~\eqref{eq:f+gh} has a solution for all but  
$$
\# \cE(\cF,\cG, \cH) =    O\(q^{1/2 + \varepsilon}\)
$$
values of $\lambda \in \F_q$, respectively.  
\end{cor}

\begin{cor}
\label{cor:EFGH2}
For any fixed 
$\varepsilon>0$, there exists $\delta > 0$ such that for  any  
sets $ \cF,\cG, \cH  \subseteq  \F_q$ which
satisfy
$$
  \# \cF  \ge q^{1/2 + \varepsilon} \qquad \text{and}  \qquad 
\# \cG \# \cH \ge  q^{2-\delta}
$$  
the 
equation~\eqref{eq:f+gh} has a solution for all but  
$$
\# \cE(\cF,\cG, \cH) =    O\(q^{\varepsilon}\)
$$
values of $\lambda \in \F_q$, respectively.  
\end{cor}

As far as we are aware, before the present work, 
multiplicative character sums have not be
used in questions of this kind. 

%

\section{Proof of Theorem~\ref{thm:2 prod}}


\subsection{Estimating $N(\cA, \cB, \cC, \cD;\lambda)$}

Let us denote by $\cX$ the set of all multiplicative characters of 
 $\F_q^*$ and by  $\cX^*$ the set of all nontrivial   characters.
For the trivial character $\chi_0$ we define $\chi_0(0) =1$ 
and put $\chi(0) = 0$ for all other characters  $\chi \in \cX^*$.

Using the orthogonality property of characters, see~\cite[Equation~(5.4)]{LN},
we write  for  the number of solutions  $N(\cA, \cB, \cC, \cD;\lambda)$ to
the equation~\eqref{eq:ab+cd}
\begin{eqnarray*}
N(\cA, \cB, \cC, \cD;\lambda) & =& \frac{1}{q-1} \sum_{a\in \cA}\sum_{b\in \cB}\sum_{c\in \cC}\sum_{d\in \cD}
\sum_{\chi\in \cX}\chi(ab - \lambda) \overline{\chi}(cd)\\
& = & \frac{\# \cA  \# \cB \# \cC \#  \cD }{q-1}+
\frac{1}{q-1} \sum_{\chi\in \cX^*}\\
& & \qquad \qquad  \sum_{a\in \cA}\sum_{b\in \cB} \chi(ab - \lambda) 
\sum_{c\in \cC}\overline{\chi}(c)\sum_{d\in \cD}
\overline{\chi}(d). 
\end{eqnarray*}
Let 
$$
W = \max_{\chi\in \cX^*} \left| \sum_{a\in \cA}\sum_{b\in \cB} \chi(ab - \lambda) \right|. 
$$
Then, using the Cauchy inequality (and again the orthogonality property of characters),
we obtain
\begin{eqnarray*}
\lefteqn{  \left|\sum_{\chi\in \cX^*} \sum_{a\in \cA}\sum_{b\in \cB} \chi(ab - \lambda) 
\sum_{c\in \cC}\overline{\chi}(c)\sum_{d\in \cD} \overline{\chi}(d)\right|} \\
& & \qquad \le W   \sum_{\chi\in \cX^*}  \left|
\sum_{c\in \cC}\overline{\chi}(c)\right|\left|\sum_{d\in \cD} \overline{\chi}(d) \right|\\
& &\qquad < W  \( \sum_{\chi\in \cX}  \left|
\sum_{c\in \cC}\overline{\chi}(c)\right|^2\)^{1/2} 
\( \sum_{\chi\in \cX} \left|\sum_{d\in \cD} \overline{\chi}(d) \right|\)^{1/2}\\
& &\qquad =  W   (q-1) \(\# \cC \#  \cD\)^{1/2}.
\end{eqnarray*}
Therefore 
\begin{equation}
\label{eq:N and W}
\left|N(\cA, \cB, \cC, \cD;\lambda)-  \frac{\# \cA  \# \cB \# \cC \#  \cD }{q-1}\right| 
<    W \(\#\cC \#  \cD\)^{1/2}.
\end{equation}

We now recall a result of  A.~A.~Karatsuba, 
see~\cite{Kar1} or~\cite[Chapter~VIII, Problem~9]{Kar2}, 
(which in turn follows from the
Weil bound and the H{\"o}lder inequality) asserting that for any integer $r \ge 1$, we
have
\begin{equation}
\label{eq:K-bound r}
W = O\( (\# \cA)^{1-1/2r} \#\cB p^{1/4r} +  
(\# \cA)^{1-1/2r} (\#\cB)^{1/2} p^{1/2r}\) ,
\end{equation}
where the implied constant depends only on $r$. 
In particular, taking $r =\rf{\varepsilon^{-1}}$,  
we see that for any $\varepsilon>0$ there exists $\delta >0$ 
such that under the condition~\eqref{eq:A and B} we have
\begin{equation}
\label{eq:K-bound}
W \le \# \cA  \#\cB q^{- 2\delta} 
\end{equation}
provided that $q$ is large enough. 
Substituting~\eqref{eq:K-bound} in~\eqref{eq:N and W} and recalling~\eqref{eq:C and D}, we obtain
the desired estimate for $N(\cA, \cB, \cC, \cD;\lambda)$.


\subsection{Estimating $T(\cA, \cB, \cC, \cD)$}

We proceed 
exactly the same way, 
however, instead of $W$, the bound on the number of solutions depends
on 
$$
V = \max_{\chi\in \cX^*} \left| \sum_{a\in \cA}\sum_{b\in \cB} \chi(a+b) \right|,  
$$
for which the same bound as~\eqref{eq:K-bound} holds as well.
 
\section{Concluding Remarks}

\subsection{Links with the Sum-Product
Problem}

We have already mentioned that  the equation~\eqref{eq:a+b+cd}
appears in the argument of M.~Garaev~\cite{Gar2} on the sum-product
problem in  finite fields. We present this argument in 
a slightly more general form. 
For two sets $\cX, \cY \subseteq \F_q$ we consider their 
sum and product sets
$$
\cU = \{x+y: \ : x \in \cX, y \in \cY \}\qquad 
\text{and}\qquad
\cV = \{xy: \ : x \in \cX, y \in \cY \}.
$$
The argument of~\cite{Gar2}
(presented there in the special case $\cX = \cY$)
is based on the observation that the equation
$$
v x_1^{-1} + x_2 = u, \qquad x_1,x_2 \in \cX, \ u \in \cU ,
\ v \in \cV, 
$$
has at least $(\# \cX)^2 \#\cY$ solutions of the form
$(x_1,x_2,u,v) = (x_1,x_2,x_2+y,x_1y)$. Combining this 
observation with~\eqref{eq:Sark T} (and assuming that 
$\#\cY \ge \# \cX$)
 we obtain
$$
\# \cU \#\cV \ge C_0 \min\left\{p\max\{\# \cX, \#\cY\}, 
\frac{(\# \cX)^2(\# \cY)^2}{p}\right\}, 
$$
for some absolute constant $C_0 > 0$, 
which for $\cX = \cY$ coincides with the result of 
M.~Garaev~\cite[Theorem~1]{Gar2}. 
It would be interesting to see whether our approach to 
estimating the number of solution to~\eqref{eq:a+b+cd}
allows to obtain stronger estimates.

\subsection{Possible Improvements}

Certainly for every concrete value of $\varepsilon$ 
one can use the bound~\eqref{eq:K-bound r}, 
instead of its simplified form~\eqref{eq:K-bound},
and optimised the choice of $r$.

Clearly~\eqref{eq:HaIo 2} is not implied by~\eqref{eq:A and B} 
and~\eqref{eq:C and D}. On the other hand, we show that our approach also
gives  an alternative prove of the corresponding 
results of~\cite{HaIo} and~\cite{Sark}
for the equations~\eqref{eq:ab+cd} and~\eqref{eq:a+b+cd}, respectively. 
In turn, using that the well-known the bound
$$
 W \le    \(q\# \cA \#  \cB\)^{1/2}, 
$$
see~\cite[Exercise~8.c]{Vin} 
we easily derive~\eqref{eq:HaIo N} from~\eqref{eq:N and W}. 
The bound~\eqref{eq:Sark T}  follows the same way from the
inequality
$$
V \le    \(q\# \cA \#  \cB\)^{1/2}.
$$ 
Certainly, if more information about the arithmetic structure 
of the sets $\cA$ and $\cB$ is known then better bounds 
can be available, see for example, estimates which are due to 
J.~Friedlander and H.~Iwaniec~\cite{FI}
and to A.~A.~Karatsuba~\cite{Kar1}. 
Moreover, in the case where $q=p$ is prime and 
$\cA = \cB = \cC  = \cB = \{1, \ldots, H\}$, 
a slight modification of the proof of~\cite[Theorem~12]{AhmShp}
shows that~\eqref{eq:ab+cd} has $H^4/p + O\(H^2 p^{o(1)}\)$
solutions for every $\lambda \in \F_p^*$,
This is nontrivial starting with $H \ge p^{1/2 + \varepsilon}$
for any $\varepsilon > 0$ and sufficiently large $p$. 
Several more results of the similar flavour 
are given by M.~Z.~Garaev and V.~Garcia~\cite{GarGar}.

\subsection{Further Problems}

Unfortunately, our approach does not seem to extend to 
the general equation~\eqref{eq:sum ab}. More precisely, combining 
bounds of exponential and multiplicative character sums, some 
results can be obtained, but they are weaker than 
those of~\cite{HaIo,HaIoKoRu}.
So, it would certainly be interesting to find a way of using bounds
of multiplicative character sums to obtain
an alternative
proof of the result of~\cite{HaIo} on the solvability 
of the equation~\eqref{eq:sum ab} under the condition~\eqref{eq:HaIo n}. 
Such a proof is likely to lead to the solvability of this
equation under several more conditions as well.

We also remark that the equation~\eqref{eq:ab+cd} (after a sign change and renaming 
some variables)  transforms into the determinant equation
$$
\det\(\begin{matrix} a&b\\c&d\end{matrix}\) =\lambda, \qquad a\in\cA,\ b \in \cB, \ c\in \cC, \ d \in \cD. 
$$
This suggests also to consider  higher dimensional determinant equations
$$
\det\(a_{ij}\)_{i,j=1}^n =\lambda, \qquad a_{ij}\in\cA_{ij}, \ i,j=1, \ldots, n,
$$
for $n^2$ sets $\cA_{ij}\subseteq \F_q$, $i,j=1, \ldots, n$. 
For structural sets, such as intervals, this and similar
questions have been studied in~\cite{AhmShp}, however its methods
do not seem to apply to the case of arbitrary sets. 

An analogue of the bound~\eqref{eq:K-bound r} has been 
given in~\cite{Shp} for character sum over points of 
an elliptic curve over $\F_q$, which has also been applied to
the studying an elliptic curve analogue of the equation~\eqref{eq:a+b+cd}.
One can also consider some mixed cases involving points on
an  elliptic curve   and field elements.

\end{document}